\documentclass[12pt]{article}
\usepackage{amsfonts,amsmath,amscd}
\usepackage[a4paper,left=3cm,right=2cm,top=3cm,bottom=2cm]{geometry}
\usepackage{hyperref}
\usepackage{lmodern}
\usepackage[utf8x]{inputenx} %
\usepackage[LGR, T1]{fontenc}
\usepackage{textcomp} 
\usepackage{graphicx}
\usepackage{wrapfig}

\newcommand{\veps}{\varepsilon}

\newcommand{\R}{\mathbb{R}}

\newtheorem{teorema}{Theorem}[section]
\newtheorem{coro}{Corollary}[section]

\newcommand{\n}{\noindent}

\begin{document}

\title{On Warped product gradient Ricci Soliton }

\author{
\textbf{M\'arcio Lemes de Sousa }
\\
{\small\it ICET - CUA, Universidade Federal de Mato Grosso,}\\
{\small\it Av. Universit\'aria nº 3.500, Pontal do Araguaia, MT,
Brazil }
\\
{\small\it e-mail:  marciolemesew@yahoo.com.br}
  \\
\textbf{Romildo Pina \footnote{ Partially supported by
CAPES-PROCAD.}}
\\
{\small\it IME, Universidade Federal de Goi\'as,}\\
{\small\it Caixa Postal 131, 74001-970, Goi\^ania, GO, Brazil }\\
{\small\it e-mail: romildo@ufg.br }
}

\maketitle

\thispagestyle{empty}

\markboth{abstract}{abstract}
\addcontentsline{toc}{chapter}{abstract}

\begin{abstract}
\noindent
In this paper we consider $M = B\times_{f}F$ warped product gradient Ricci solitons.
We show that when the Hessian matrix of $f$ is not zero then the potential function depends only on the base and the fiber $F$ is necessarily Einstein manifold. Moreover,  we provide all such solutions in the case of steady gradient Ricci solitons when the base is conformal to an $n$-dimensional pseudo-Euclidean space, invariant under the action of an $(n-1)$-dimensional translation group, and the fiber $F$ is Ricci-flat.
\end{abstract}


\section{Introduction and main statements}

Let $(M, g)$ be a semi-Riemannian manifold of dimension $n\geq 3$. We say that $(M, g)$
is a \emph{ \bf gradient Ricci soliton} if there exists a differentiable  function $h:M\rightarrow\mathbb{R}$
(called the potential function) such that
\[
\mbox{Ric}_g+\mbox{Hess}_g(h)=\rho g, \qquad \rho\in\R,
\]
where Ric$_g$ is the Ricci tensor, Hess$_g(h)$ is the Hessian of $h$
with respect to the metric $g$, and $\rho$ is a real number.
A gradient Ricci soliton is said to be {\em shrinking, steady or expanding} if
$\rho>0$, $\rho=0$ or $\rho<0$, respectively. When $M$ is a Riemannian
manifold usually one requires the manifold to be complete. In the case of
semi-Riemannian manifold, one does not require $(M,g)$ to be complete.

Ricci solitons are natural extensions of Einstein manifolds and they appear
as self-similar solutions of the Ricci flow $\partial g(t)/\partial t =-2\mbox{Ric}g(t)$.
Moreover, Ricci solitons are important in understanding singularities of the Ricci flow.

Simple examples of gradient Ricci solitons are obtained by
considering $\R^n$ with the canonical metric $g$. Then
$(\R^n,g)$  is a gradient Ricci soliton, where  $h(x)=A|x|^2/2+g(x,B)+C$,
$A,C\in \R$ and $B\in \R^n$ are all the potential functions.
In this case, $(\R^n,g)$ is a shrinking, steady or expanding soliton according to
the sign of the constant $A$.  Bryant \cite{BRYANT} proved that there exists a
complete, steady,  gradient Ricci soliton spherically symmetric for any
$n\geq 3$, which is known as Bryant's soliton.

Recently,  Cao and Chen \cite{CAOCHEN} showed that any complete, steady, gradient Ricci
soliton, locally conformally flat, up to homothety, is either isometric to
the Bryant's soliton or is flat. Complete, shrinking gradient
solitons, conformally flat, have been characterized as being quotients of
$\R^n$, $S^n$ or $\R\times S^{n-1}$ (see \cite{FG}).

Although the Ricci soliton equation was introduced and studied initially in the Riemannian context,
Lorentzian Ricci solitons have been recently investigated in \cite{BBGG}, \cite{BCGG}
and \cite{Onda}, where the authors show that there are important differences with the
Riemannian case. The existence of Lorentzian, steady, gradient Ricci
solitons which are locally conformally flat and distinct from Bryant's
solitons, was proved  in \cite{BBGG}. In \cite{BGG}, the authors gave a local
characterization of the Lorentzian gradient Ricci solitons which are locally
conformally flat.

In \cite{Bar}, Barbosa-Pina-Tenenblat, considered gradient Ricci solitons, conformal to an $n$-dimensional
pseudo-Euclidean space, which are invariant under the action of an $(n-1)$-dimensional translation group.
The one provided all such solutions in the case of steady gradient Ricci solitons.

Our purpose is to generalize the results in \cite{Bar}.  To get these generalizations, we
have to use warped product manifolds to study
gradient Ricci solitons that are non locally conformally flat. Then considering  $(B, g_{B})$ and $(F, g_{F})$
semi-Riemannian manifolds, and let $f>0$ be a smooth function
on $B$, the warped product $M = B\times_{f}F$ is the product
manifold  $ B\times F$ furnished with metric tensor
$$\tilde{g} = g_{B} + f^{2}g_{F},$$
$B$ is called the base of $M = B\times_{f}F$, $F$ the fiber and
$f$ is the warping function. For example, polar coordinates
determine a warped product in the case of constant curvature
spaces, the case corresponds to $\R^{+}\times_{r}S^{n - 1} $.

There are several studies correlating warped product manifolds and
locally conformally flat manifolds, see \cite{BGV1}, \cite{BGV3} and their references.

Many of the original examples of gradient Ricci solitons arise as warped products over a one dimensional base
(cf. \cite{CCGG}, \cite{KLCL}).

In this paper, initially we prove that if a warped product  $M = B\times_{f}F$ is a gradient Ricci soliton  with the   hypothesis that there is at least one pair of vector $ ( X_i, X_k) $ of the base, such that $ \mbox{Hess}_{g_{B}}(f) ( X_i, X_k) \ne 0$
 then the potential function depends only on the base and the fiber is necessarily an Einstein manifold (see Theorem 1.1 and Corollary 1.1).  The Theorem 1.1 and Corollary 1.1 generalize the results obtained in \cite{KBH}  where  the authors  studied  warped product  gradient Ricci solitons with one-dimensional base.

  In what follows,
 we consider warped product $M = B\times_{f}F$ gradient Ricci solitons, where the base is conformal to a pseudo-Euclidean space
which are invariant under the action of an $(n-1)$-dimensional translation group and the fiber is an Einstein manifold. More precisely,  let  $(\R^n,g)$ be  the pseudo-Euclidean space, $n\geq 3$, with
coordinates $x=(x_1,\cdots, x_n)$ and $g_{ij}=\delta_{ij}\veps_i$ and let $M = (\mathbb{R}^{n}, \overline{g})\times _{f}F^{m}$ be a warped
product where $\displaystyle \overline{g} =
\frac{1}{\varphi^{2}}g$,  $F$ a semi--Riemannian
Einstein manifold  with constant Ricci curvature $\lambda_{F}$, $m\geq 1$,
$f,\varphi, h:\mathbb{R}^{n}\rightarrow \mathbb{R}$,  smooth
functions, and $f$ is a positive function. In Theorem 1.2 we find necessary and sufficient
conditions for the warped product metric $\widetilde{g} = \overline{g} +
f^{2}g_{F} $ be a gradient Ricci soliton, namely
\begin{equation}\label{grad1}
\mbox{Ric}_{\tilde{g}}+\mbox{Hess}_{\tilde{g}}(h)=\rho \tilde{g}, \qquad \rho\in\R.
\end{equation}

 In Theorem 1.3, we consider $f$, $\varphi$ and $h$  invariant under the action of an (n-1)--dimensional
translation group and  let $\xi=\sum_{i=1}^{n}\alpha_ix_i, \; \alpha_i\in
 \R$, be a basic invariant for an $(n-1)$-dimensional translation group.
We want to obtain differentiable functions $\varphi(\xi), f(\xi)$ and $h(\xi)$,  such that the metric
$\widetilde{g}$ is a gradient Ricci soliton. We first obtain necessary and sufficient conditions on $f(\xi)$, $\varphi(\xi)$ and $h(\xi)$ for the existence of $\widetilde{g}$. We show
that these conditions are different depending on the direction $\alpha=\sum_{i=1}^n \alpha_i\partial/\partial x_i$ being null (lightlike) or not. We observe that in the null case the metrics $\widetilde{g}$ and $g_{F}$
are necessarily gradient Ricci solitons steady and Ricci-flat, respectively.

Considering  $M =  (\mathbb{R}^{n}, \overline{g})\times _{f}F^{m}$ and $F$ Ricci-flat manifold,  we obtain all the metrics $\widetilde{g} = \overline{g} + f^{2}g_{F} $, which are gradient Ricci solitons steady and are invariant under the action of an $(n-1)$-dimensional translation group.
We prove that if the direction $\alpha$ is timelike or spacelike,  the functions $f$, $\varphi$ and $h$ depend on the dimensions $n$, $m$ and also on a finite number of parameters. In fact, the solutions are explicitly given in Theorems 1.4, 1.5 and 1.6.
 If the direction $\alpha$ is null, then there are infinitely many solutions. In fact, in this case,
for any given positive differentiable functions $f(\xi)$ and $\varphi(\xi)$,  the function $h(\xi)$ satisfies a linear ordinary differential equation of second order (see Theorem 1.7).  We illustrate this fact with some explicit examples.

When the dimension of the fiber $F$ is $m = 1$ we consider $M =
(\R^n, \overline{g})\times _{f}\mathbb{R}$ and in this case
$\lambda_{F} = 0$.

In what follows, we state our main results.
We denote by $\varphi_{ ,x_ix_j}$, $f_{ ,x_ix_j}$ and $h_{ ,x_ix_j}$
the second order derivative of $\varphi$, $f$ and $h$,  with respect to
$x_ix_j$.

\begin{teorema}\label{teor1}
Let  $M = B^{n}\times_{f}F^{m}$  be a warped product semi-Riemannian manifold with metric $\widetilde{g}$.  If
the warped product  $\widetilde{g}= g_{B} + f^{2}g_{F}$  is a gradient Ricci soliton with    $h:M \rightarrow\mathbb{R}$ as potential function, and there is at least one pair of vector $ ( X_i, X_k) $ of the base, such that $ \mbox{Hess}_{g_{B}}(f) ( X_i, X_k) \ne 0$, 
then  $h$  depends only on the base. 

\end{teorema}

\begin{coro}
Let $M = B\times_{f}F$ be a warped product semi-Riemannian manifold with metric $\widetilde{g}$. If
the warped product  metric $\widetilde{g}= g_{B} + f^{2}g_{F}$ is a gradient Ricci soliton with  $h:B \rightarrow\mathbb{R}$
as potential function, and $f$ is non-constant, then the fiber is an Einstein manifold.
\end{coro}
The results obtained in Theorem 1.1 and Corollary 1.1 are valid in the  Riemannian case.
Motivated by the previous results we study the problem considering warped product with $h$ depending
only on the base and the fiber an Einstein manifold.

\begin{teorema}\label{teor2}
Let $(\R^n,g)$ be a pseudo-Euclidean space, $n\geq 3$ with
coordinates $x=(x_1,\cdots, x_n)$ and $g_{ij}=\delta_{ij}\veps_i$.
Consider $M = (\mathbb{R}^{n}, \overline{g})\times _{f}F^{m}$ a
warped product, where $\displaystyle \overline{g} =
\frac{1}{\varphi^{2}}g$,  $F$ a semi--Riemannian Einstein manifold
with constant Ricci curvature $\lambda_{F}$, $m\geq 1$,
$f,\varphi, h:\mathbb{R}^{n}\rightarrow \mathbb{R}$,  smooth
functions and $f$ is a positive function. Then  the warped
product metric $\widetilde{g} = \overline{g} + f^{2}g_{F} $ is a gradient Ricci soliton with $h$
as potential function if, and only if, the functions $f$, $\varphi$ and $h$ satisfy:

\begin{equation}\label{eqphij}
(n - 2)f\varphi_{ ,x_{i}x_{j}} +f\varphi h_{ ,x_{i}x_{j}}  - m\varphi f_{ ,x_{i}x_{j}} -
m\varphi_{ ,x_{i}}f_{ ,x_{j}} - m\varphi_{ ,x_{j}}f_{ ,x_{i}} + f\varphi_{ ,x_{i}}h_{ ,x_{j}} +
f\varphi_{ ,x_{j}}h_{ ,x_{i}} = 0,
\end{equation}
where $1\leq i\neq j\leq n,$

\begin{equation}\label{eqphii}
\begin{array}{rll}
\varphi\big[(n-2)f\varphi_{ ,x_ix_i} + f\varphi h_{ ,x_{i}x_{i}}- m\varphi f_{ ,x_{i}x_{i}} -
2m\varphi_{ ,x_{i}}f_{ ,x_{i}}+ 2f\varphi_{ ,x_{i}}h_{ ,x_{i}}\big] &+&
\\\varepsilon_{i}\displaystyle\sum_{k
=1}^{n}\varepsilon_{k}\big[f\varphi\varphi_{ ,x_{k}x_{k}} - (n
-1)f\varphi_{,x_{k}}^{2}
  + m\varphi \varphi_{ ,x_{k}}f_{ ,x_{k}} -f\varphi
  \varphi_{ ,x_{k}}h_{ ,x_{k}}\big]
  &=& \varepsilon_{i}\rho f, 1\leq i \leq n
\end{array}
\end{equation}
and
\begin{equation}\label{eqphll}
\sum_{k =1}^{n}\varepsilon_{k}\big[-f\varphi^{2} f_{ ,x_{k}x_{k}} + (n-2)f\varphi f_{ ,x_{k}}\varphi_{,x_{k}} - (m-1)\varphi^{2} f_{,x_{k}}^{2} + f\varphi^{2} f_{,x_{k}}h_{,x_{k}}\big]
 = \rho f^{2} - \lambda_{F}.
\end{equation}
\end{teorema}

We want to find solutions of the system (\ref{eqphij}),
(\ref{eqphii}) and (\ref{eqphll})
 of the form $\varphi(\xi)$, $f(\xi)$ and $h(\xi)$, where   $\xi=\sum_{i=1}^{n}\alpha_ix_i, \; \alpha_i\in
 \R$. Whenever $\sum_{i=1}^{n}\veps_i\alpha_i^2\neq 0$, without loss of generality, we may consider $\sum_{i=1}^{n}\veps_i\alpha_i^2=\pm 1$.
The following theorem provides the system of ordinary differential
equations that must be satisfied by such solutions.

\vspace{.2in}

\noindent{\bf Theorem 1.3.} {\em Let $( \R^n, g)$ be a pseudo-Euclidean space, $n\geq 3$, with coordinates $x=(x_1,\cdots, x_n)$ and $g_{ij}=\delta_{ij}\veps_i$.
 Consider $M = (\mathbb{R}^{n}, \overline{g})\times_{f}F^{m}$,
 where $\displaystyle \overline{g} = \frac{1}{\varphi^{2}}g$, $F^{m}$ a semi--Riemannian  Einstein manifold
 with constant Ricci curvature $\lambda_{F}$ and smooth
functions $\varphi(\xi)$,  $f(\xi)$ and $h(\xi)$, where
$\xi=\sum_{i=1}^{n}\alpha_ix_i, \; \alpha_i\in
 \R,$ and  $\sum_{i=1}^{n}\veps_i\alpha_i^2=\veps_{i_0}$ or $\sum_{i=1}^{n}\veps_i\alpha_i^2=0$. Then  the warped product metric $\widetilde{g} = \overline{g} + f^{2}g_{F} $ is a gradient Ricci soliton with $h$
as potential function if, and only if, the functions $f$, $\varphi$ and $h$ satisfy:\\ }
 \\
{\em  (i)
 \begin{equation} \label{solitonxi}
\left\{ \begin{array}{lcl}
 (n-2)f\varphi'' + f\varphi h'' - m\varphi f''- 2m\varphi' f' + 2f\varphi' h' = 0, \\
\displaystyle\sum_{k =
1}^{n}\varepsilon_{k}\alpha_{k}^{2}\big[f\varphi\varphi''
- (n-1)f\varphi'^{2} + m\varphi\varphi'f' - f\varphi\varphi'h'\big] =  \rho f\\
\displaystyle\sum_{k =
1}^{n}\varepsilon_{k}\alpha_{k}^{2}\big[-f\varphi^{2}f'' +
(n-2)f\varphi\varphi'f' - (m-1)\varphi^{2}f'^{2} + f\varphi^{2}f'h'\big] = \rho f^{2}
- \lambda_{F},
\end{array} \right.
 \end{equation}
  whenever }$\displaystyle\sum_{i=1}^{n}\veps_i\alpha_i^2=\veps_{i_0}, $ and  \\
{\em (ii)
\begin{equation}\label{eqxi0}
 (n-2)f\varphi''+ f\varphi h''- m\varphi f''- 2m\varphi' f' + 2f\varphi'h'=0, \quad \mbox{ \em and }\quad  \rho = \lambda_{F}=0,
 \end{equation}
  whenever } $\displaystyle\sum_{i=1}^{n}\veps_i\alpha_i^2=0. $

In the following two results we provide all the solutions of
(\ref{solitonxi}) when $\widetilde{g}$ is steady gradient Ricci solitons, i.e., $\rho = 0$, and
 $F^{m}$ is a Ricci-flat manifold, i.e., $\lambda_{F} = 0$.

 \vspace{.2in}
 \noindent{\bf Theorem 1.4.} {\em Let $( \R^n, g)$ be a
pseudo-Euclidean space, $n\geq 3$, with coordinates
$x=(x_1,\cdots, x_n)$ and $g_{ij}=\delta_{ij}\veps_i$.
 Consider  smooth functions $\varphi(\xi)$,  $f(\xi)$ and $h(\xi)$, where
$\xi=\sum_{i=1}^{n}\alpha_ix_i, \; \alpha_i\in
 \R,$  $\sum_{i=1}^{n}\veps_i\alpha_i^2=\veps_{i_0}$, given by
 \begin{equation} \label{steady_xi_N}
 \left\{\begin{array}{lcl}
 \varphi_{\pm}(\xi)&=& \displaystyle\frac{c_{2}}{(N_{\pm}\xi + b)^{\frac{k}{N_{\pm}}}}, \\
 \mbox{} \\
 f_{\pm}(\xi)&=& \displaystyle\frac{c_{1}}{(N_{\pm}\xi + b)^{\frac{1}{N_{\pm}}}},\\
 \mbox{}\\
 h_{\pm}(\xi)& =& \displaystyle-\frac{m-(n-2)k + N_{\pm}}{N_{\pm}}\ln\left( N_{\pm}\xi + b \right)
 \end{array} \right.
\end{equation}
where $k, c_{1}, c_{2}, N_{\pm}$ and $b$ are constant with $k, c_{1}, c_{2}>0$ and $N_{\pm} = -k \pm\sqrt{m + k^{2}(n - 1)}$.
 Then the warped product $M = (\R^n, \overline{g})\times _{f}F^{m}$, with $\overline{g} = \frac{1}{\varphi^{2}}g$ and $F$ Ricci-flat,
 is a steady gradient Ricci soliton with $h$ as potential function. These solutions are defined on the half space determined by
$\sum_{i=1}^{n}\alpha_ix_i< \frac{b}{-k + \sqrt{m + k^{2}(n - 1)}}$ or $\sum_{i=1}^{n}\alpha_ix_i> \frac{b}{-k - \sqrt{m + k^{2}(n - 1)}}\cdot$}

\vskip10pt
\noindent{\bf Theorem 1.5} {\em Let $(
 \R^n, g)$ be a pseudo-Euclidean space,  $n\geq 3$, with
coordinates $x=(x_1,\cdots, x_n)$ and $g_{ij}=\delta_{ij}\veps_i$.
 Consider smooth functions $x(\xi)$, and $z(\xi)$, where $\xi= \sum_{i=1}^{n}\alpha_ix_i, \; \alpha_i\in \R$,  $\sum_{i=1}^{n}\veps_i\alpha_i^2=\pm 1$, given by }
 \begin{equation} \label{steady_xi_W}
  x = c_{3}\displaystyle\left(z + k -\sqrt{m + k^{2}(n - 1)}\right)^{\frac{a - 1}{2}}\left(z + k +\sqrt{m + k^{2}(n - 1)}\right)^{-\frac{a + 1}{2}}
 \end{equation}
 \begin{equation}\label{steady xi z}
 z'= \displaystyle -c_{3}\displaystyle\left(z + k -\sqrt{m + k^{2}(n - 1)}\right)^{\frac{a + 1}{2}}\left(z + k +\sqrt{m + k^{2}(n - 1)}\right)^{-\frac{a - 1}{2}}
 \end{equation}
{\em where $a, k, c_{3}$ are constants with, $k,
c_{3}>0$ and $a = \displaystyle \frac{k}{\sqrt{m + k^{2}(n -1)}}$.\\ Let  $\varphi(\xi)$,  $f(\xi)$ and $h(\xi)$ be functions obtained by integrating,
\begin{equation}\label{so1}
	\frac{\varphi'}{\varphi}(\xi) = k\frac{f'}{f}(\xi),\ \ \ \   \frac{f'}{f}(\xi) = x(\xi),\ \ \ \ h'(\xi) = [z(\xi) + m -k(n-2)]x(\xi)\cdot
\end{equation}
 Then the warped product $M = (\R^n, \overline{g})\times _{f}F^{m}$, with $\overline{g} = \frac{1}{\varphi^{2}}g$ and $F$ Ricci-flat, is a steady gradient Ricci soliton with $h$ as potential function.}

We remark that the obtained  metric, in the Theorems 1.4 and 1.5, are non locally conformally flat.

 \noindent{\bf Theorem 1.6} {\em Let $(
 \R^n, g)$ be a pseudo-Euclidean space,  $n\geq 3$, with
coordinates $x=(x_1,\cdots, x_n)$ and $g_{ij}=\delta_{ij}\veps_i$. Consider $M = (\mathbb{R}^{n}, \overline{g})\times _{f}F^{m}$ a warped product where $\displaystyle \overline{g} =
\frac{1}{\varphi^{2}}g$. Let $\widetilde{g} = \overline{g} + f^{2}g_{F} $, be a steady gradient Ricci soliton with $h$ as potential function and $F$ Ricci-flat. Then $\varphi$, $f$ and $h$ are invariant under an (n-1)-dimensional translation group whose basic invariant is $\xi = \sum_{i =1}^{n}\alpha_{i}x_{i}$, where $\alpha = \sum_{i = 1}^{n}\alpha_{i}\partial/\partial x_{i}$ is a non null vector if, and only if, $\varphi$, $f$ and $h$ are given as in Theorems $1.4$ or $1.5$. }

 The following theorem shows that there are infinitely many warped products $M = (\R^n, \overline{g})\times _{f}F^{m}$
steady gradient Ricci solitons with $h$ as potential function, which are invariant under the action of an
$(n-1)$--dimensional group acting on $\mathbb{R}^{n}$, when $\alpha
=  \sum_{i = 1}^{n}\alpha_{i}\partial/\partial\
x_{i}$ is null vector.

\noindent{\bf Theorem 1.7} {\em Let $\varphi(\xi) $ and $f(\xi)$ be any positive
differentiable functions, where $\xi =  \sum_{i =
1}^{n}\alpha_{i}x_{i}$ and $ \sum_{i =
1}^{n}\varepsilon_{i}\alpha_{i}^{2} = 0$. Then the
function $h(\xi)$ given by
\begin{equation}\label{so2}
h(\xi) = \int\left(\varphi^{-2}\left[\int(m\frac{f''}{f}\varphi^{2} + 2m\varphi\varphi'\frac{f'}{f} - (n-2)\varphi\varphi'')d\xi + c_{4}\right]\right)d\xi + c_{5}; c_{4}, c_{5}\in\R
\end{equation}

satisfies  (\ref{eqxi0}) and $M = (\R^n,
\overline{g})\times _{f}F^{m}$ is a steady gradient Ricci soliton with $h$ as potential function. }

Before proving our main results, we present two examples illustrating
Theorem 1.7. Let $f(\xi) = \varphi(\xi) = ke^{A\xi}$, where $k>0$, $A\neq 0$,
$\xi =  \sum_{i = 1}^{n}\alpha_{i}x_{i}$ and
$ \sum_{i = 1}^{n}\varepsilon_{i}\alpha_{i}^{2} = 0$.
 Solving  (\ref{so2}) we get
$$h(\xi) = (3m - (n-2))\frac{A}{2}\xi - \frac{c_{4}}{2Ak^{2}}e^{-2A\xi} + c_{5};\ \  c_{4}, c_{5}\in\mathbb{R}.$$ It follows from Theorem 1.7 that $M = (\R^n,
\overline{g})\times _{f}F^{m}$ is a steady gradient Ricci soliton, with $h$ as potential function. Similarly, if we choose $f(\xi) = e^{-\xi^{2}}$ and $\varphi(\xi) = e^{\xi}$, then it follows from Theorem 1.7 that
$$h(\xi) = \frac{2m}{3}\xi^{3} - 2m\xi^{2} + \left(m - \frac{n-2}{2}\right)\xi - \frac{1}{2}c_{4}e^{-2\xi} + c_{5}; \ \ c_{4}, c_{5}\in\mathbb{R},$$
is the potential function of the steady gradient Ricci soliton $\tilde{g} = \bar{g} + f^{2}g_{F}$.

 \section{Proofs of the Main Results}

 \n \textbf{Proof of Theorem \ref{teor1}:}
   
   Let  $M = B\times_{f}F$  be a gradient Ricci soliton with potential $h$, then 
  \[
 \mbox{Ric}_{\tilde{g}}+\mbox{Hess}_{\tilde{g}}(h)=\rho \tilde{g}, \qquad \rho\in\R.
 \]
  Considering  $ X_{1},\ldots X_{n}\in{\cal L}(B)$  $Y_{1}\ldots, Y_{m}\in{\cal L}(F)$, where  ${\cal L}(B)$ and ${\cal L}(F)$ are respectively the lift of a vector field on $B$ and $F$ to $B\times F$, we have
\begin{eqnarray}
\label{Ric1}
\left\{\begin{array}{lll}
  Ric_{\widetilde{g}}(X_{i},X_{j}) &=& Ric_{g_{B}}(X_{i},X_{j}) - \frac{m}{f}Hess_{g_{B}}f(X_{i},X_{j}),\ \forall\ i,\ j = 1,\ldots n \\
    Ric_{\widetilde{g}}(X_{i},Y_{j}) &=& 0, \forall\ i= 1,\ldots n,\ j = 1,\ldots m  \\
    Ric_{\widetilde{g}}(Y_{i},Y_{j}) &=& Ric_{g_{F}}(Y_{i}, Y_{j}) -
  \widetilde{g}(Y_{i}, Y_{j})\left(\frac{\triangle_{g_{B}}f}{f} + (m-1)\frac{\widetilde{g}(\nabla f, \nabla
  f)}{f^{2}}\right),\ \forall\  i,\ j = 1,\ldots m
\end{array}
\right.
\end{eqnarray}

  How  $\tilde{g}(X_{i}, Y_{j}) = 0$ for $i = 1\ldots n$  and $j = 1\ldots m$ (see \cite{O'neil}), we have that
 \begin{equation}\label{Hess}
 \mbox{Hess}_{\tilde{g}}(h)(X_{i}, Y_{j})= 0.
 \end{equation}
 Note that
 \[\mbox{Hess}_{\tilde{g}}(h)(X_{i}, Y_{j}) = X_{i}Y_{j}(h) - (\nabla_{X_{i}}Y_{j})(h), i= 1\ldots n, j = 1\ldots m, \]

  where   $\nabla$ is the connection of   $M$. Since (see  \cite{O'neil}),
 \[\nabla_{X_{i}}Y_{j} = \frac{X_{i}(f)}{f}Y_{j}\]
  we have 
 \begin{equation}\label{Hess2}
 \mbox{Hess}_{\tilde{g}}(h)(X_{i}, Y_{j}) = h,_{x_{i}y_{j}} -\frac{f,_{x_{i}}}{f}h,_{y_{j}}.
 \end{equation}

  Using  (\ref{Hess2}) and (\ref{Hess}) we obtain:
 \begin{equation}\label{Hes}
  h,_{x_{i}y_{j}} -\frac{f,_{x_{i}}}{f}h,_{y_{j}} = 0.
 \end{equation}
  If   $h = h(x_{1},\ldots,x_{n})$  then the equation  (\ref{Hes}) is trivially satisfied.   Suppose that there is at least one $y_j$,  with  $ 1 \leq  j \leq m $   such that  $ h,_{y_j} \ne 0$.
   In this case, it follows from (\ref{Hes}) that
  \begin{equation}\label{Hess3}
  \frac{h,_{x_{i}y_{j}}}{h,_{y_{j}}} = \frac{f,_{x_{i}}}{f}, \forall\ \ i = 1\ldots n,\ \ j = 1\ldots m.
  \end{equation}
integrating (\ref{Hess3}) in relation to   $x_{i}$we obtain:
\[\ln{h,_{y_{j}}} = \ln{f} + l(\hat{x_{i}})\]
this is,
\begin{equation}\label{Hess4}
h,_{y_{j}} = fe^{l(\hat{x_{i}})}.
\end{equation}
 Fixing $i$ and  $j$ in  ( \ref{Hess4}) and deriving in relation to  $x_{k}$ with $k\neq i$, we obtain 
 \[h,_{y_{j}x_{k}} = f,_{x_{k}}e^{l(\hat{x_{i}})} +f l,_{x_{k}}e^{l(\hat{x_{i}})}\]
 
 which is equivalent to
 
 \[\frac{h_{x_{k}y_{j}}}{h_{y_{j}}} = \frac{f_{x_{k}}}{f} + l_{x_{k}},\]
 and again using   (\ref{Hess3}) we have:
 \[l_{x_{k}} = 0,\]
 This means that  $l$ does not depend $x_{k}$, this is 
 \[l = l(\hat{x_{i}}, \hat{x_{k}}).\]
Repeating this process  we obtain  that $l$  depends only on the fiber, i.e.
 \[l = l(y_{1},\ldots, y_{m}),\]
 therefore,
 \begin{equation}\label{Hess5}
 h_{y_{j}} = fe^{l(y_{1},\ldots, y_{m})}  \quad \quad   with \quad  \quad  1 \leq  j \leq m   .
 \end{equation}

 Integrating  (\ref{Hess5}) in relation to $y_{j}$,  we obtain:
 \begin{equation}\label{Hess6}
 h(x_{1},\ldots,x_{n};y_{1}\ldots y_{m}) = f(x_{1}, \ldots, x_{n})\int e^{l(y_{1},\ldots, y_{m})}dy_{j} + m(\hat{y_{j}}).
 \end{equation}

 Using the first equation  (\ref{Ric1}), we have,
 \[ \mbox{Hess}_{\tilde{g}}(h)(X_{i}, X_{k})= \rho g_{B}(X_{i},X_{k}) - Ric_{g_{B}}(X_{i},X_{k}) + \frac{m}{f}Hess_{g_{B}}f(X_{i},X_{k}),\ \forall\ i,\ k = 1,\ldots n, \]
proving that $\mbox{Hess}_{\tilde{g}}(h)(X_{i}, X_{k})$,  $ \forall  i,\ k = 1,\ldots n$,   depends only on the base.  Thus,  considering  $j$ fixed  in equation  (\ref{Hess6})  we have:
\begin{equation}\label{Hess7}
\displaystyle\frac{\partial}{\partial y_{j}}\mbox{Hess}_{\tilde{g}}(h)(X_{i}, X_{k}) = 0, \forall\ i,\ k = 1,\ldots n.
\end{equation}
On the other hand  $  \forall  i,\ k = 1,\ldots n$  we have
\[
\frac{\partial}{\partial y_{j}}\mbox{Hess}_{\tilde{g}}(h)(X_{i}, X_{k}) = \frac{\partial}{\partial y_{j}}\left(\mbox{Hess}_{\tilde{g}}(f)(X_{i},X_{k})\int e^{l(y_{1},\ldots, y_{m})}dy_{j} + \mbox{Hess}_{\tilde{g}}(m)(X_{i},X_{k})\right)\]
Since  $m = m(\hat{y_{j}})$,    and using the definition of the Hessian get that   

\begin{equation}\label{Hess32}
\displaystyle\frac{\partial}{\partial y_{j}}\mbox{Hess}_{\tilde{g}}(m)(X_{i}, X_{k}) = 0, \forall\ i,\ k = 1,\ldots n, 
\end{equation}

and as  $ \mbox{Hess}_{\tilde{g}}(f)(X_{i},X_{k}) =  \mbox{Hess}_{g_{B}}(f)(X_{i},X_{k}) $  follow that 
 
\begin{equation}\label{Hess8}
\frac{\partial}{\partial y_{j}}\mbox{Hess}_{\tilde{g}}(h)(X_{i}, X_{k}) = \mbox{Hess}_{g_{B}}(f)(X_{i},X_{k}) e^{l(y_{1},\ldots, y_{m})}, \forall\ i,\ k = 1,\ldots n.
\end{equation}

Using  (\ref{Hess7}) and (\ref{Hess8}) we obtain that:
\[\mbox{Hess}_{g_{B}}(f)(X_{i},X_{k}) e^{l(y_{1},\ldots, y_{m})} = 0, \forall\ i,\ k = 1,\ldots n,\]

We have for  hypothesis that there is at least one pair of vector $ ( X_i, X_k) $ of the base, such that $\mbox{Hess}_{g_{B}}(f) ( X_i, X_k) \ne 0$ then

 \[e^{l(y_{1},\ldots, y_{m})} = 0,\]
 
 but this is impossible. Therefore     $ h,_{y_{j}} = 0   \quad \quad   \forall  \quad  \quad  j = 1\ldots m$. Consequently  $h$ depends only on the base.

  This concludes the proof  of the Theorem \ref{teor1}.

\hfill $\Box$

\vspace{.2in}

 \n \textbf{Proof of Corollary 1.1}
Let $M = B\times_{f}F$ be a gradient Ricci soliton with potential $h$ defined only on the base, then by the equation (\ref{grad1}) we have,
  \[
 \mbox{Ric}_{\tilde{g}}+\mbox{Hess}_{\tilde{g}}(h)=\rho \tilde{g}, \qquad \rho\in\R.
 \]
Considering  $Y, Z\in{\cal L}(F)$, we have that
\begin{equation}\label{s10}
 \mbox{Ric}_{\tilde{g}}(Y, Z)+\mbox{Hess}_{\tilde{g}}(h)(Y, Z)=\rho \tilde{g}(Y, Z),
  \end{equation}
 but, \[\tilde{g}(Y, Z) = f^{2}g_{F}(Y, Z)\] and \[\mbox{Ric}_{\tilde{g}}(Y, Z) = \mbox{Ric}_{g_{F}}(Y, Z) - (f\triangle_{g_{B}}f + (m -1 )\|grad_{g_{B}}f\|^{2})g_{F}(Y, Z)\]
   (see for example \cite{O'neil}). Replacing $\mbox{Ric}_{\tilde{g}}(Y, Z)$ in (\ref{s10}), we have
 \begin{equation}\label{hess2}
 \mbox{Ric}_{g_{F}}(Y, Z) = (f^{2} + f\triangle_{g_{B}}f + (m -1 )\|grad_{g_{B}}f\|^{2})g_{F}(Y, Z) - \mbox{Hess}_{\tilde{g}}(h)(Y, Z).
 \end{equation}
It follows from (\ref{hess2}) that $F$ is Einstein if, and only if,
 \[\mbox{Hess}_{\tilde{g}}(h)(Y, Z) = \lambda g_{F}(Y, Z).\]
Indeed, as by hypotheses $h$ depends only on the base, thus $grad_{\tilde{g}}h = grad_{g_{B}}h$,  then
 \[\mathcal{H}(grad_{\tilde{g}}h) =  grad_{g_{B}}h, \]
 \[\mathcal{V}(grad_{\tilde{g}}h) = 0\]
 and
 \[\nabla_{Y}(grad_{\tilde{g}}h) = \frac{grad_{g_{B}}h(f)}{f}Y.\]
 Therefore
\begin{equation}\label{hessf}
\mbox{Hess}_{\tilde{g}}(h)(Y, Z) =  \frac{grad_{g_{B}}h(f)}{f}\tilde{g}(Y, Z)= f grad_{g_{B}}h(f)g_{F}(Y, Z).
\end{equation}

This concludes the proof Corollary 1.1.

\hfill $\Box$

\vspace{.2in}

\n \textbf{Proof of Theorem \ref{teor2}:} Assume initially that $m>1$. It follows from
\cite{O'neil} that if $X_{1}, X_{2},\ldots, X_{n}\in{\cal
L}(\mathbb{R}^{n})$ and $Y_{1}, Y_{2}, \ldots, Y_{m}\in{\cal
L}(F)$ (${\cal L}(\mathbb{R}^{n})$ and ${\cal L}(F)$ are respectively the lift of a vector field on $\mathbb{R}^{n}$ and $F$ to $\mathbb{R}^{n}\times F$ ), then
\begin{eqnarray}
\label{ric1}
\left\{\begin{array}{lcl}
Ric_{\widetilde{g}}(X_{i},X_{j}) &=& Ric_{\overline{g}}(X_{i},X_{j}) - \frac{m}{f}Hess_{\overline{g}}f(X_{i},X_{j}),\ \forall\ i,\ j = 1,\ldots n \\
    Ric_{\widetilde{g}}(X_{i},Y_{j}) &=& 0, \forall\ i= 1,\ldots n,\ j = 1,\ldots m  \\
    Ric_{\widetilde{g}}(Y_{i},Y_{j}) &=& Ric_{g_{F}}(Y_{i}, Y_{j}) -
  \widetilde{g}(Y_{i}, Y_{j})(\frac{\triangle_{\overline{g}}f}{f} + (m-1)\frac{\widetilde{g}(\nabla f, \nabla
  f)}{f^{2}}),\ \forall\  i,\ j = 1,\ldots m
\end{array}
\right.
\end{eqnarray}

 It is well known (see, e.g., \cite{Be}) that if $\overline{g}=\frac{1}{\varphi^2}g$  , then
\[
Ric_{\overline{g}}=\frac{1}{\varphi^2}\left \{(n-2)\varphi
Hess_{g}(\varphi)+[\varphi \Delta_g
\varphi-(n-1)|\nabla_g\varphi|^2]g \right \}\,.
\]
Since $g(X_{i},X_{j}) = \varepsilon_{i}\delta_{ij}$, we have

\begin{eqnarray*}
  Ric_{\overline{g}}(X_{i}, X_{j}) &=& \frac{1}{\varphi}\left \{(n-2)
Hess_{g}(\varphi)(X_{i}, X_{j}) \right \}\ \forall\ i\neq\ j = 1,\ldots n  \\
  Ric_{\overline{g}}(X_{i}, X_{i}) &=& \frac{1}{\varphi^2}\left \{(n-2)\varphi
Hess_{g}(\varphi)(X_{i}, X_{i})+[\varphi \Delta_g
\varphi-(n-1)|\nabla_g\varphi|^2]\varepsilon_{i} \right \}\ \forall\
i=1,\ldots n.
\end{eqnarray*}
 Since $Hess_{g}(\varphi)(X_{i}, X_{j}) = \varphi_{,x_{i}x_{j}}$
, $\Delta_g \varphi = \displaystyle\sum_{k=1}^{n}\varepsilon_{k}\varphi_{,x_{k}x_{k}}$ and
$|\nabla_g\varphi|^2 = \displaystyle\sum_{k =
1}^{n}\varepsilon_{k}\varphi_{,x_{k}}^{2} $, we have
\begin{equation}\label{ric2}
\left\{ \begin{array}{ccl}
  Ric_{\overline{g}}(X_{i}, X_{j}) &=& \displaystyle\frac{(n-2)\varphi_{,x_{i}x_{j}}}{\varphi}\qquad \forall\ i\neq\ j:1\ldots\ n \\
  Ric_{\overline{g}}(X_{i}, X_{i}) &=& \displaystyle\frac{(n-2)\varphi_{,x_{i}x_{i}} +
  \varepsilon_{i}\displaystyle\sum_{k=1}^{n}\varepsilon_{k}\varphi_{,x_{k}x_{k}}}{\varphi}
  - (n - 1)\varepsilon_{i}\sum_{k=1}^{n}\frac{\varepsilon_{k}\varphi_{,x_{k}}^{2}}{\varphi^{2}}
\end{array}
\right.
\end{equation}

Recall that
\[
Hess_{\overline{g}}(f)(X_{i},X_{j})=f_{ ,x_ix_j}-\sum_k
\overline{\Gamma}_{ij}^k f_{,x_k},
\]
where $\overline{\Gamma}_{ij}^k$ are the Christoffel symbols of the metric $\overline{g}$. For $i,\ j,\ k$ distinct, we have
\[
\overline{\Gamma}_{ij}^k= 0\ \ \ \ \ \ \ \ \ \overline{\Gamma}_{ij}^i= -\frac{\varphi_{,x_{j}}}{\varphi}\ \ \ \ \ \ \ \ \ \ \overline{\Gamma}_{ii}^k= \varepsilon_{i}\varepsilon_{k}\frac{\varphi_{,x_{k}}}{\varphi}\ \ \ \ \ \ \ \ \ \overline{\Gamma}_{ii}^i= -\frac{\varphi_{,x_{j}}}{\varphi} \] therefore,

\begin{equation}\label{hes1}
\left\{
\begin{array}{ccc}
Hess_{\overline{g}}(f)(X_{i},X_{j}) &=& \displaystyle
f_{,x_ix_j}+\frac{\varphi_{,x_j}}{\varphi}f_{,x_i}
+\frac{\varphi_{,x_i}}{\varphi}f_{,x_j},\ \forall\ i\neq j = 1\ldots n\\
Hess_{\overline{g}}(f)(X_{i},X_{i}) &=& \displaystyle
f_{,x_ix_i}+2\frac{\varphi_{,x_i}}{\varphi}f_{,x_i} -
\varepsilon_{i}\sum_{k=1}^{n}\varepsilon_{k}\frac{\varphi_{,x_k}}{\varphi}f_{,x_k}
\end{array}
\right.
\end{equation}

 Substituting (\ref{ric2}) and (\ref{hes1}) in the first equation of the (\ref{ric1}) we  obtain
 \begin{equation}\label{ric3}
  Ric_{\widetilde{g}}(X_{i},X_{j}) =
 \frac{(n-2)\varphi_{,x_{i}x_{j}}}{\varphi}-\frac{m}{f}\left[f_{,x_ix_j}+\frac{\varphi_{,x_j}}{\varphi}f_{,x_i}
+\frac{\varphi_{,x_i}}{\varphi}f_{,x_j}\right],\ \forall\ i\neq j
\end{equation}

and

\begin{eqnarray}\label{ric4}
  Ric_{\widetilde{g}}(X_{i},X_{i}) &=& \frac{(n-2)\varphi_{,x_{i}x_{i}} +
  \varepsilon_{i}\displaystyle\sum_{k=1}^{n}\varepsilon_{k}\varphi_{,x_{k}x_{k}}}{\varphi}
  - (n -
  1)\varepsilon_{i}\sum_{k=1}^{n}\frac{\varepsilon_{k}\varphi_{,x_{k}}^{2}}{\varphi^{2}} \nonumber\\
  &-&  \frac{m}{f}\left[f_{,x_ix_i}+2\frac{\varphi_{,x_i}}{\varphi}f_{,x_i} -
\varepsilon_{i}\sum_{k=1}^{n}\varepsilon_{k}\frac{\varphi_{,x_k}}{\varphi}f_{,x_k}.\right]
\end{eqnarray}

On the other hand,
\begin{equation}\label{gr}
\left\{
\begin{array}{ccl}
  Ric_{{g_{F}}}(Y_{i}, Y_{j}) &=& \lambda_{F}g_{F}(Y_{i}, Y_{j}) \\
  \widetilde{g}(Y_{i}, Y_{j}) &=& f^{2}g_{F}(Y_{i}, Y_{j}) \\
  \Delta_{\overline{g}}f &=& \varphi^{2}\sum_{k = 1}^{n}\varepsilon_{k}f_{,x_{k}x_{k}} - (n-2)\varphi\sum_{k = 1}^{n}\varepsilon_{k}\varphi_{,x_k}f_{,x_k} \\
  \widetilde{g}(\nabla f, \nabla f) &=& \varphi^{2}\sum_{k =
  1}^{n}\varepsilon_{k}f_{,x_{k}}^{2}
\end{array}
\right.
\end{equation}

Substituting (\ref{gr}) in the third equation of (\ref{ric1}),
we have

\begin{equation}\label{ric5}
 Ric_{\widetilde{g}}(Y_{i},Y_{j}) = \gamma_{ij}g_{F}(Y_{i}, Y_{j})
\end{equation}
where, \[\gamma_{ij } =\lambda_{F} -
f\varphi^{2}\sum_{k = 1}^{n}\varepsilon_{k}f_{,x_{k}x_{k}} +
(n-2)f\varphi\sum_{k = 1}^{n}\varepsilon_{k}
\varphi_{,x_k}f_{,x_k} - (m - 1)\varphi^{2}\sum_{k =
1}^{n}\varepsilon_{k}f_{,x_{k}}^{2}  .\]

We want to find $h$ satisfying
\begin{equation}\label{rsg}
\mbox{Ric}_{\tilde{g}}(X, Y)+\mbox{Hess}_{\tilde{g}}(h)(X, Y)=\rho \tilde{g}(X, Y), \forall X, Y\in \mathcal{X}(M).
\end{equation}

On the other hand, since $h:\mathbb{R}^{n}\rightarrow \mathbb{R}$, we have that

\[Hess_{\widetilde{g}}(h)(X_{i}, X_{j}) = Hess_{\overline{g}}(h)(X_{i}, X_{j}), \forall 1\leq i, j\leq n\]
i.e.,
\begin{equation}\label{hess3}
\left\{
\begin{array}{ccc}
Hess_{\widetilde{g}}(h)(X_{i},X_{j}) &=& \displaystyle
h_{,x_ix_j}+\frac{\varphi_{,x_j}}{\varphi}h_{,x_i}
+\frac{\varphi_{,x_i}}{\varphi}h_{,x_j},\ \forall\ i\neq j = 1\ldots n\\
Hess_{\widetilde{g}}(h)(X_{i},X_{i}) &=& \displaystyle
h_{,x_ix_i}+2\frac{\varphi_{,x_i}}{\varphi}h_{,x_i} -
\varepsilon_{i}\sum_{k=1}^{n}\varepsilon_{k}\frac{\varphi_{,x_k}}{\varphi}h_{,x_k}.
\end{array}
\right.
\end{equation}
By substituting (\ref{ric3}) and the first equation of (\ref{hess3}) into (\ref{rsg}), we obtain (\ref{eqphij}). Again replacing (\ref{ric4}) and the second equation of (\ref{hess3}) into (\ref{rsg}) we get (\ref{eqphii}).
Now for $X_{i}\in{\cal L}(\mathbb{R}^{n})$ and $Y_{j}\in{\cal L}(F)$ ($1\leq i \leq n$ and $1\leq j\leq m$) we get
\[Hess_{\widetilde{g}}(h)(X_{i}, Y_{j}) = 0\]
by Theorem \ref{teor1}. In this case the equation (\ref{rsg}) is trivially satisfied.
Taking $Y_{i}, Y_{j}\in{\cal L}(F)$ with $1\leq i,j\leq m$ by using the equation (\ref{hessf}) we have:
\[\mbox{Hess}_{\tilde{g}}(h)(Y_{i}, Y_{j}) =  \frac{grad_{\overline{g}}h(f)}{f}\tilde{g}(Y_{i},Y_{j})= f grad_{\overline{g}}h(f)g_{F}(Y_{i}, Y_{j})\]
but,
\[ grad_{\overline{g}}h(f) = \overline{g}(grad_{g_{B}}h, grad_{g_{B}}f) = \varphi^{2}\sum_{k=1}^{n}\varepsilon_{k}f_{,x_{k}}h_{,x_{k}}\]
Then we have,
\begin{equation}\label{hess4}
\mbox{Hess}_{\tilde{g}}(h)(Y_{i}, Y_{j}) = f\varphi^{2}\sum_{k=1}^{n}\varepsilon_{k}f_{,x_{k}}h_{,x_{k}}g_{F}(Y_{i}, Y_{j}).
\end{equation}

By substituting (\ref{ric5}) and (\ref{hess4}) into (\ref{rsg}) we obtain (\ref{eqphll})

The reciprocal of this theorem can be easily verified.
 In the case  $m = 1$ just remember that

\begin{eqnarray*}
  Ric_{\widetilde{g}}(X_{i},X_{j}) &=& Ric_{\overline{g}}(X_{i},X_{j}) - \frac{1}{f}Hess_{\overline{g}}f(X_{i},X_{j}),\ \forall\  i,\ j = 1,\ldots n \\
  Ric_{\widetilde{g}}(X_{i},Y) &=& 0,\ \forall\ i= 1,\ldots n\\
  Ric_{\widetilde{g}}(Y, Y) &=& - \widetilde{g}(Y, Y)\frac{\triangle_{\overline{g}}f}{f}.
\end{eqnarray*}
In this case the equation (\ref{eqphij}) and (\ref{eqphii}) remain the same and the equation (\ref{eqphll}) reduces to

\[ -\varphi^{2}\sum_{k =1}^{n}\varepsilon_{k} f_{ ,x_{k}x_{k}} + (n
-2)\varphi\sum_{k =1}^{n}\varepsilon_{k} f_{ ,x_{k}}\varphi_{
,x_{k}} + \varphi^{2}\sum_{k=1}^{n}\varepsilon_{k}f_{,x_{k}}h_{,x_{k}} = \rho f .\]

This concludes the proof Theorem \ref{teor2}.

\hfill $\Box$

\vspace{.2in}

\n \textbf{Proof of Theorem 1.3} Assume initially that $m>1$. Let  $\overline{g}=\varphi^{-2}g$ be a conformal metric of $g$. We are assuming that $\varphi(\xi)$, $f(\xi)$ and $h(\xi)$ are functions of $\xi$, where $\xi=\sum \limits_{i=1}^n\alpha_ix_i$, $\alpha_i\in \R$ and $\sum_i\veps_i\alpha_i^2=\veps_{i_0}$ or $\sum_i\veps_i\alpha_i^2=0$. Hence, we have
\[
\varphi_{,x_i}=\varphi'\alpha_i, \qquad
\varphi_{,x_ix_j}=\varphi''\alpha_i\alpha_j,
\]
\[
f_{,x_i}=h'\alpha_i, \qquad f_{,x_ix_j}=f''\alpha_i\alpha_j,
\]
and
\[
h_{,x_i}=h'\alpha_i, \qquad h_{,x_ix_j}=h''\alpha_i\alpha_j\,.
\]
Substituting  these expressions into (\ref{eqphij}), we get
\[
(n - 2)f\varphi''\alpha_i\alpha_j + f\varphi h''\alpha_i\alpha_j - m\varphi
f''\alpha_{i}\alpha_{j} - 2m\varphi'f'\alpha_{i}\alpha_{j} + 2f\varphi'h'= 0,
\qquad \forall\ i\neq j.
\]

If there exist $i\neq j$ such that $\alpha_i\alpha_j\neq 0$, then this equation reduces to
\begin{equation}  \label{eqphiLL}
(n - 2)f\varphi'' + f\varphi h'' - m\varphi f'' - 2m\varphi'f' + 2f\varphi'h'= 0
\end{equation}
Similarly, considering  equation (\ref{eqphii}), we get
\begin{equation}\label{st1}
\begin{array}{rll}
 \alpha_i^2 \varphi\left[ (n-2)f\varphi'' + f\varphi h'' - m\varphi f'' - 2m\varphi'f' + 2f\varphi'h' \right]
 &+& \\\veps_i\sum_k \veps_k\alpha_k^2 \left[f\varphi \varphi'' -(n-1)f(\varphi')^2 + m\varphi \varphi'f' - f\varphi\varphi'h'  \right]=\veps_i\rho f.
 \end{array}
\end{equation}
Due to the relation between $\varphi''$, $f''$ and $h''$ given in (\ref{eqphiLL}), the equation (\ref{st1}) reduces to
\begin{equation}\label{inter}
\sum_k \veps_k\alpha_k^2 \left[f\varphi \varphi''
-(n-1)f(\varphi')^2 + m\varphi \varphi'f' - f\varphi\varphi'h' \right]=\rho f.
\end{equation}

 Analogously to equation (\ref{eqphll}) reduces to
 \begin{equation}
\sum_k \veps_k\alpha_k^2 \left[-f\varphi^{2}f''
+(n-2)f\varphi\varphi' f' - (m - 1)\varphi^{2} f'^{2} + f\varphi^{2} f'h' \right] = \rho.
f^{2} - \lambda_{F}
\end{equation}

Therefore, if $\sum_k \veps_k\alpha_k^2=\veps_{i_0}$, we obtain the equations of the system (\ref{solitonxi}). If $\sum_k
\veps_k\alpha_k^2=0$, we have (\ref{eqphiLL})  satisfied and
(\ref{inter}) implies $\rho=0$, hence $\lambda_{F} = 0$ i.e.,
(\ref{eqxi0}) holds.

If for all $i\neq\ j$, we have $\alpha_{i}\alpha_{j} = 0$, then $\xi = x_{i_{0}}$, and equation (\ref{eqphij}) is trivially satisfied for all $i\neq\ j$. Considering (\ref{eqphii}) for $i\neq\ i_{0}$, we get
\[\displaystyle\sum_{k =
1}^{n}\varepsilon_{k}\alpha_{k}^{2}[f\varphi\varphi''
- (n-1)f\varphi'^{2} + m\varphi\varphi'f' - f\varphi\varphi'h'] = \rho f\]
and hence, the second equation of (\ref{solitonxi}) is satisfied. Considering $i = i_{0}$ in (\ref{eqphii}) we get that the first equation of (\ref{solitonxi})  is satisfied.

Considering $i = i_{0}$ or $i \neq\ i_{0}$ in (\ref{eqphll}), we get that the third equation (\ref{solitonxi}) is  satisfied.

When $m = 1$  the first and the second equation of the system (\ref{solitonxi}) are the same, and the third equation  reduces to

\[\displaystyle\sum_{k =
1}^{n}\varepsilon_{k}\alpha_{k}^{2}[-\varphi^{2}f'' +
(n-2)\varphi\varphi'f' +\varphi^{2}f'h']  = \rho f.\]

This concludes the proof Theorem 1.3. \label{}

\hfill $\Box$

\vspace{.2in}

In order to prove Theorems 1.4 and 1.5, we consider functions $\varphi(\xi), f(\xi)$ and $h(\xi)$, where  $\xi \displaystyle = \sum_{i=1}^{n}\alpha_{i}x_{i}, \alpha_{i}\in\mathbb{R}, \sum_{i=1}^{n}\varepsilon_{i}\alpha_{i}^{2} = \pm 1$. It follows from Theorem 1.3 that $\widetilde{g} = \overline{g} + f^{2}g_{F}$ is a steady gradient Ricci soliton with Ricci-flat $F$ and $h$ as potential function if, and only if, $\varphi, f$ and $h$ satisfy
\begin{equation}
\left\{ \begin{array}{lcl}
\displaystyle (n-2)\frac{\varphi''}{\varphi} + h'' - m\frac{f''}{f} - 2m\frac{\varphi'}{\varphi} \frac{f'}{f} + 2\frac{\varphi'}{\varphi} h' = 0, \\
\displaystyle \frac{\varphi''}{\varphi} - (n-1)\left(\frac{\varphi'}{\varphi}\right)^{2} + m\frac{\varphi'}{\varphi}\frac{f'}{f} - \frac{\varphi'}{\varphi}h' =  0\\
\displaystyle-\frac{f''}{f} + (n-2)\frac{\varphi'}{\varphi}\frac{f'}{f} - (m-1)\left(\frac{f'}{f}\right)^{2} + \frac{f'}{f}h' = 0,
\end{array} \right.
 \end{equation}
 i.e.,
 \begin{equation}\label{st2}
\left\{ \begin{array}{lll}
\displaystyle (n-2)\left(\frac{\varphi'}{\varphi}\right)' +  (n-2)\left(\frac{\varphi'}{\varphi}\right)^{2}+ h'' - m\left(\frac{f'}{f}\right)' - m\left(\frac{f'}{f}\right)^{2} - 2m\frac{\varphi'}{\varphi} \frac{f'}{f} + 2\frac{\varphi'}{\varphi} h' = 0, \\
\displaystyle \left(\frac{\varphi'}{\varphi}\right)' - (n-2)\left(\frac{\varphi'}{\varphi}\right)^{2} + m\frac{\varphi'}{\varphi}\frac{f'}{f} - \frac{\varphi'}{\varphi}h' =  0\\
\displaystyle -\left(\frac{f'}{f}\right)'- m\left(\frac{f'}{f}\right)^{2} + (n-2)\frac{\varphi'}{\varphi}\frac{f'}{f}  + \frac{f'}{f}h' = 0.
\end{array} \right.
 \end{equation}
Then, it follows from the second equation of (\ref{st2}) that
\begin{equation}\label{po1}
h' =\displaystyle \frac{(\frac{\varphi'}{\varphi})'}{\frac{\varphi'}{\varphi}} - (n-2)\frac{\varphi'}{\varphi} + m\frac{f'}{f},
\end{equation}
on the other hand, by using the third equation of (\ref{st2}) we have
\begin{equation}\label{po2}
h' = \displaystyle\frac{(\frac{f'}{f})'}{\frac{f'}{f}} - (n-2)\frac{\varphi'}{\varphi} + m\frac{f'}{f},
\end{equation}
hence substituting equation (\ref{po2}) into (\ref{po1}), we get
\begin{equation}\label{po3}
\displaystyle \frac{(\frac{\varphi'}{\varphi})'}{\frac{\varphi'}{\varphi}} = \frac{(\frac{f'}{f})'}{\frac{f'}{f}}.
\end{equation}
Integrating the equation (\ref{po3}) we obtain:
\begin{equation}\label{po4}
\frac{\varphi'}{\varphi} = k\frac{f'}{f}, \ \ k>0.
\end{equation}
Substituting (\ref{po4}) in the system (\ref{st2}), we get

\begin{equation}\label{st3}
\left\{ \begin{array}{lll}
\displaystyle h'' +[(n-2)k - m]\left(\frac{f'}{f}\right)' = [2mk + m - (n-2)k^{2}] \left(\frac{f'}{f}\right)^{2} - 2k \frac{f'}{f}h' \\
\displaystyle\left(\frac{f'}{f}\right)' = [(n-2)k - m]\left(\frac{f'}{f}\right)^{2} + \frac{f'}{f}h'\cdot
\end{array} \right.
 \end{equation}
Considering
\begin{equation}\label{mu}
x(\xi) = (f'/f)(\xi) \qquad and \qquad y =  h' +[(n-2)k - m]\frac{f'}{f}
\end{equation}
the system of equations (\ref{st3}) is equivalent to
\begin{equation}\label{st4}
\left\{ \begin{array}{lll}
 y' = [m + (n-2)k^{2}]x^{2} - 2k xy \\
x' = xy\;\cdot
\end{array} \right.
 \end{equation}
It follows from (\ref{st4}) that
\begin{equation}\label{st5}
xydy - \{[m + (n-2)k^{2}]x^{2} - 2kxy\}dx = 0\cdot
\end{equation}
In order to study this equation (see e.g. \cite{HT} p. 37), we take
\begin{equation}\label{st6}
y(\xi) = x(\xi)z(\xi)\cdot
\end{equation}
Theorem 1.4 and 1.5 are obtained by considering $z$ to be a non zero constant and a nonconstant function, respectively.

\noindent {\bf Proof of Theorem 1.4.} We consider solutions of the system (\ref{st4}), as in (\ref{st6}) where $z(\xi) = N$ and $N$ is a nonzero constant, i.e., $y(\xi) = Nx(\xi)$. By substituting $y$ in the first and second equations of (\ref{st4}), we get
 \[ x' =\displaystyle \frac{m + (n-2)k^{2} -2kN}{N}x^{2}\]
and
\[x' = Nx^{2}\]
respectively. Comparing the two expressions we conclude that $N$ must satisfy
\[N^{2} +2kN - (m + (n-2)k^{2}) = 0\cdot\]
 Therefore, we get two values of $N$, given by
 \begin{equation}
 N_{\pm} = -k \pm\sqrt{m + (n-1)k^{2}}\cdot
 \end{equation}
 Going back to the second equation, we have that
\[\displaystyle \frac{x'}{x^{2}} = N_{\pm}\cdot\]
 Hence we have
 \[x_{\pm}(\xi)\displaystyle =\frac{-1}{N_{\pm}\xi + b}\qquad and \qquad y_{\pm}(\xi) = \frac{-N_{\pm}}{N_{\pm}\xi + b}. \]
From (\ref{po4}) and (\ref{mu}) we have $\displaystyle \frac{\varphi'}{\varphi} = k\frac{f'}{f}$, $x(\xi) = (f'/f)(\xi)$ and $y =\displaystyle  h' +[(n-2)k - m]\frac{f'}{f}$, hence we get $\varphi, f$ and $h$ given by (\ref{steady_xi_N}).

This concludes the proof Theorem 1.4.

\hfill $\Box$

\noindent {\bf Proof of Theorem 1.5.}
 We now consider solutions of the system (\ref{st4}), where $y(\xi) = x(\xi)z(\xi)$ and $z(\xi)$ is a smooth nonconstant function. Substituting $y$ into (\ref{st5}) and assuming, without loss of generality, $x\neq 0$ on an open subset, we get
 \begin{equation}\label{st7}
 \displaystyle \frac{dx}{x} = - \frac{z}{z^{2} + 2kz - m - (n-2)k^{2}}dz \cdot
 \end{equation}
 Integrating this equation, we get $x$ in terms of $z$ as given by equation (\ref{steady_xi_W}). Since $y = xz$, it follows from the first equation of (\ref{st4}) that
 \[x'z + xz' = [m + (n-2)k^{2}]x^{2} - 2k x^{2}z,\]
 and second equation of (\ref{st4}) that
 \[x' = x^{2}z\cdot\]
Hence we get,
\begin{equation}\label{st8}
z' = x[m + (n-2)k^{2} -2kz - z^{2}].
\end{equation}
 Substituting the equation (\ref{steady_xi_W}) into (\ref{st8}), we conclude that $z$ must satisfy the differential equation (\ref{steady xi z}). Now a straight forward computation shows that both equations of (\ref{st4}) are satisfied. Since we have determined $z(\xi)$ and $x(\xi)$, we can obtain $\varphi(\xi), f(\xi)$ and $h(\xi)$ integrating (\ref{so1}).

This concludes the proof Theorem 1.5.

 \hfill $\Box$

 \noindent {\bf Proof of Theorem 1.6.}
 When $\rho = \lambda_{F} = 0$, by introducing the auxiliary functions $\displaystyle\frac{\varphi'}{\varphi}(\xi) = k\frac{f'}{f}(\xi),   \frac{f'}{f}(\xi) = x(\xi)$ and $h'(\xi) = [z(\xi) + m -k(n-2)]x(\xi)$, we have seen that (\ref{solitonxi}) is equivalent to the system (\ref{st4}) for $x$ and $y$. The solutions of this system can be written as $y(\xi) = x(\xi)z(\xi)$, where $z(\xi)$ is a nonzero function.

 If $z(\xi)$ is a non zero constant, then the proof of Theorem 1.4 shows that the solutions of (\ref{st4}) are given by (\ref{steady_xi_N}). If the function is not constant, then proof of Theorem 1.5 shows that $z$ is determined by (\ref{steady xi z}) and $x$ is given algebraically in terms of $z$ by (\ref{steady_xi_W}). Then one gets the functions $\varphi, f$ and the potential $h$ by integrating the ordinary differential equations given by (\ref{so1}). This concludes the proof Theorem 1.6.

 \hfill $\Box$

\noindent {\bf Proof of Theorem 1.7.}
Let $(\R^n,g)$ be a pseudo-Euclidean space, $n\geq 3$ with
coordinates $x=(x_1,\cdots, x_n)$ and $g_{ij}=\delta_{ij}\veps_i$. Consider $M = (\mathbb{R}^{n}, \overline{g})\times _{f}F^{m}$ a warped product where $\displaystyle \overline{g} =\frac{1}{\varphi^{2}}g$,  $F$ a semi--Riemannian manifold Ricci--flat. Let $\varphi(\xi)$ and $f(\xi)$ be any non-vanishing differentiable functions invariant by the translation of $(n-1)$-dimensional translation group, whose basic invariant is $\xi =  \displaystyle \sum_{i=1}^{n}\alpha_{i}x_{i}$, where $\alpha_{i}\in \R$ and $\displaystyle\sum_{i=1}^{n}\veps_i\alpha_i^2=0, $. Then it follows from Theorem 1.3 that the warped product metric $\widetilde{g} = \overline{g} + f^{2}g_{F}$ is a gradient Ricci soliton with $h$ as a potential if, and only if, $\rho = 0$ and $h$ satisfies the linear ordinary differential equation (\ref{eqxi0}) determined by $\varphi$ and $f$. Then is easy to see that
\[h'(\xi) =\varphi^{-2}\left[\int(m\frac{f''}{f}\varphi^{2} + 2m\varphi\varphi'\frac{f'}{f} - (n-2)\varphi\varphi'')d\xi + c_{4}\right] \]
and hence is given by (\ref{so2})

 This concludes the proof Theorem 1.7.

 \hfill $\Box$

\end{document}